\documentclass[letterpaper, 10 pt, conference]{ieeeconf}  

\IEEEoverridecommandlockouts                              
\usepackage{cite}
\usepackage{amsmath,amssymb,amsfonts}
\usepackage{algorithmic}
\usepackage{graphicx}
\usepackage{textcomp}
\usepackage{xcolor}
\usepackage{mathrsfs}  
\usepackage[font=scriptsize]{caption}
\let\proof\relax

\usepackage{amsthm,bm}
\usepackage{subfigure}
\newtheorem{theorem}{Theorem}
\newtheorem{lemma}{Lemma}

\newcommand{\cA}{\mathcal{A}}
\newcommand{\cB}{\mathcal{B}}
\newcommand{\cC}{\mathcal{C}}
\newcommand{\cQ}{\mathcal{Q}}
\newcommand{\cR}{\mathcal{R}}
\newcommand{\cJ}{\mathcal{J}}
\newcommand{\cK}{\mathcal{K}}
\newcommand{\calL}{\mathcal{L}}
\newcommand{\cL}{\mathscr{L}}

\newcommand{\cH}{\mathcal{H}}

\newcommand{\T}{\mathbb{T}}
\newcommand{\R}{\mathbb{R}}
\newcommand{\C}{\mathbb{C}}

\newcommand{\E}{\mathbb{E}}

\def\BibTeX{{\rm B\kern-.05em{\sc i\kern-.025em b}\kern-.08em
    T\kern-.1667em\lower.7ex\hbox{E}\kern-.125emX}}

\title{\LARGE \bf
Complete Decentralization of Linear Quadratic Gaussian Control for the Discrete Wave Equation}

\author{Addie McCurdy and Emily Jensen
\thanks{This work was not supported by any organization}
\thanks{A. McCurdy is with the Department of Applied Mathematics, University of Colorado Boulder, Boulder, CO 80309, USA {\tt\small addie.mccurdy@colorado.edu}}%
\thanks{E. Jensen is with Department of Electrical, Computer \& Energy Engineering, University of Colorado Boulder, Boulder, CO 80309, USA {\tt\small ejensen@colorado.edu}}%
}

\begin{document}

\maketitle
\thispagestyle{empty}
\pagestyle{empty}

\begin{abstract}

   The linear quadratic Gaussian (LQG) control problem for the linear wave equation on the unit circle with fully distributed actuation and partial state measurements is considered. An analytical solution to a spatial discretization of the problem is obtained. The main result of this work illustrates that for specific parameter values, the optimal LQG policy is \emph{completely decentralized}, meaning only a measurement at spatial location $i$ is needed to compute an optimal control signal 
   to actuate
   at this location. The relationship between performance and decentralization as a function of parameters is explored. Conditions for complete decentralization are related to metrics of kinetic and potential energy quantities and control effort.


\end{abstract}

\section{Introduction}
As sensing and actuating capabilities of engineering systems grow, it is increasingly important to leverage distributed or decentralized control policies to mitigate potential computational or communication bottlenecks. It is known to some extent that properties of the plant will inevitably impact the potential performance of a distributed policy. On one hand, the optimal feedback policy for certain spatially-invariant plants will be approximately localized \cite{bamieh2002distributed, jensen2020localization}; in contrast, certain plant structures are innately uncontrollable by certain classes of distributed controllers \cite{sezer1981structurally}.

Toward an answer to how plant structures and parameters impact the potential performance of a distributed or decentralized control policy, we examine here the structure of the optimal linear quadratic Gaussian (LQG) controller for the discretized wave equation on the unit circle, assuming distributed actuation and distributed partial state measurements. The standard wave equation and variations of it are prevalent across engineering domains, e.g. models of transmission line dynamics \cite{remoissenet2013waves}, pressure variations in a compressible gas \cite{ockendon2016waves}, elastic materials \cite{deng2017elastic} and soft-robotic systems \cite{tanaka2012mechanics}. Mass-spring systems \cite{lindberg2024h} and simplified models of vehicular platoons \cite{jovanovic2005ill} are examples of discrete systems whose dynamics are a spatial discretization of the wave equation; incorporating damping terms and nonlinearities leads to models of systems like networks of RLC circuits \cite{gander2020asymptotic} or off-shore multi-modular solar panels \cite{al2020hydrodynamic}. We focus on the linear, undamped wave equation to obtain analytical results that will provide physical insight to numerical results obtained for more complex wave-like systems, e.g., with damping or nonlinearities. Indeed, the wave equation has been shown to provide meaningful approximations of more complex discrete systems such as power networks \cite{chakrabortty2012graph}.

This work builds on results that provide conditions for decentralization of the Kalman filter for the wave equation over the real line \cite{arbelaiz2022information}. In contrast to this recent work, we restrict to a finite spatial domain and discretize over the spatial variable. 
The finite spatial domain restriction is perhaps more realistic as physical systems do not have infinite spatial extent, but it prevents us from leveraging the approach of \cite{arbelaiz2022information} that leads to explicit bounds on spatial decay rates of the Kalman filter gain. This restriction also enables a clear characterization of the relationship between locality and performance, by avoiding an infinite measure of cost that occurs with an unbounded spatial domain \cite{arbelaiz2022information}.

Our main result provides a condition on system parameters to obtain complete decentralization of both the linear quadratic regulator (LQR) and the Kalman filtering (KF) gains for the discretized wave equation with periodic boundary conditions. By appropriately selecting system parameters, we can ensure optimal control performance with minimal communication and computation requirements. Interestingly, we highlight that this complete decentralization result can only occur for certain formulations of the LQG problem. The conditions for a formulation that is amenable to decentralization are connected to metrics of potential and kinetic energy as well as to covariance of sensor noise measurements. The sensor noise covariance condition may be tied to recent results of estimation for lower-order linear PDEs \cite{arbelaiz2020distributed}. Alternatively, we can view this formulation as a more ``accurate'' description of the full PDE which accounts for the domain of the $C_0$-semigroup \cite{curtain2020introduction} generated by the unbounded operator describing the state dynamics of the system. 


The structure and main results of this paper are summarized as follows. In Section~\ref{sec:notation}, we introduce relevant notation and preliminaries. In Section~\ref{sec:problemPDE}, we state an LQG problem for the wave equation over the unit circle, which we discretize and solve in Section~\ref{sec:LQGL2}. We then show that communication is always required for optimality {in the preceding LQG formulation}. Our main results are presented in Section~\ref{sec:smartLQG}. We provide a complete characterization of the system parameters that lead to complete decentralization of the LQR and KF gains when we implement `Sobolev-like' cost functions. In Section~\ref{sec:discussion}, we illustrate that parameters can be selected to satisfy both of these conditions, leading to maximal locality of the LQG controller. The relationship between controller performance and locality is explored.

\section{Notation and Mathematical Preliminaries}\label{sec:notation}
Let $H_1$ and $H_2$ be Hilbert spaces. We denote the space of linear operators from $H_1$ to $H_2$ as $\cL(H_1,H_2)$, or $\cL(H_1)$ if $H_1=H_2$. When $H_1,H_2$ are infinite dimensional, we denote operators in $\cL(H_1,H_2)$ with capital script letters (e.g. $\cA$). When $H_1,H_2$ are finite dimensional, operators in $\cL(H_1,H_2)$ are matrices and denoted with capital print letters (e.g. $A$). The exception to this convention will be the identity operator, which we denote as $I$ regardless of the spaces it maps between.

We consider spatially distributed systems where the spatial domain is {either the unit circle $\mathbb{T}$ or $n$ equispaced points in $\T$, denoted $\Omega\in\R^n$. Equivalently, $\mathbb{T}$ can be thought of as an interval of length $h$ with periodic boundary conditions.} When the spatial domain is $\T$, we write our spatio-temporal signals with explicit dependence on the spatial coordinate $x$, e.g. $p(x,t)$. For each $t$, $p(\cdot,t)$ is an element of an infinite-dimensional Hilbert space. For our purposes, such signals will belong to one of the following: 
\begin{itemize}
    \item $L_2(\T)$, the set of square integrable functions from $\T$ to $\C$ equipped with the inner product 
    \begin{equation}
        \langle f,g\rangle_{L_2(\T)} :=\int_{x\in\T}g(x)^*f(x)\,dx
    \end{equation}
    \item $\cH_\alpha(\T)$, the first order Sobolev space of weakly differentiable functions from $\T$ to $\C$ equipped with the inner product 
    \begin{equation}
        \langle f,g\rangle_{\cH_\alpha(\T)} :=\langle f,g\rangle_{L_2(\T)}+\alpha^2\langle \partial_xf,\partial_x g\rangle_{L_2(\T)}
    \end{equation}
    for some $\alpha>0.$ 
\end{itemize}
When the spatial domain is $\Omega$, we write the state as a time-varying vector in $\R^n$, e.g. $\bm p(t)$. In this problem, for each $t$ such signals will belong to the following Hilbert space: 
\begin{itemize}
    \item $\ell_2$, the space of vectors in $\R^n$ equipped with the inner product 
    \begin{equation}
        \langle \bm{f}, \bm{g}\rangle_{\ell_2}=\bm{g}^T\bm{f}.
    \end{equation} 
\end{itemize}

Hats denote the unitary spatial Discrete Fourier Transform (DFT) of signals in $\ell_2$:
\begin{equation}
    \hat{\bm{f}}_k(t)=(F\bm{f})(k,t):=\frac{1}{\sqrt{n}}\sum_{j=1}^{n-1}\bm{f}_j(t)\exp{\left(-2\pi i\frac{k}{n}j\right)}
\end{equation}
where $k\in\{0,\ldots,n-1\}$ is the spatial frequency. We write $\hat{\bm{f}}=[\hat{\bm{f}}_0,\ldots,\hat{\bm{f}}_{n-1}]^T\in\R^n$. The DFT can be written as multiplication by a matrix $F$ so that $\hat{\bm{f}}=F\bm{f}.$ 
Since $F$ is unitary,
\begin{equation}\label{eq:planch}
    \langle \bm{f},\bm{g}\rangle_{\ell_2}=\langle \hat{\bm{f}},\hat{\bm{g}}\rangle_{\ell_2}.
\end{equation}

$\mathrm{diag}(f(k))$ denotes a diagonal matrix with $k^{th}$ diagonal element given by $f(k)$ where $f$ is any function of $k$.

A \textit{circulant matrix} is completely determined by its first row, as each sequential row is the previous row shifted right one index. The DFT diagonalizes circulant matrices, i.e. if $M$ is an $n\times n$ circulant matrix, then \[\hat{M}=FMF^{-1}=\mathrm{diag}(\lambda_M(k))\] where $\lambda_M(k), k=1,\ldots, n$ is the $k^{th}$ eigenvalue of $M$. 

\begin{figure}
    \centering
    \includegraphics[width=\linewidth]{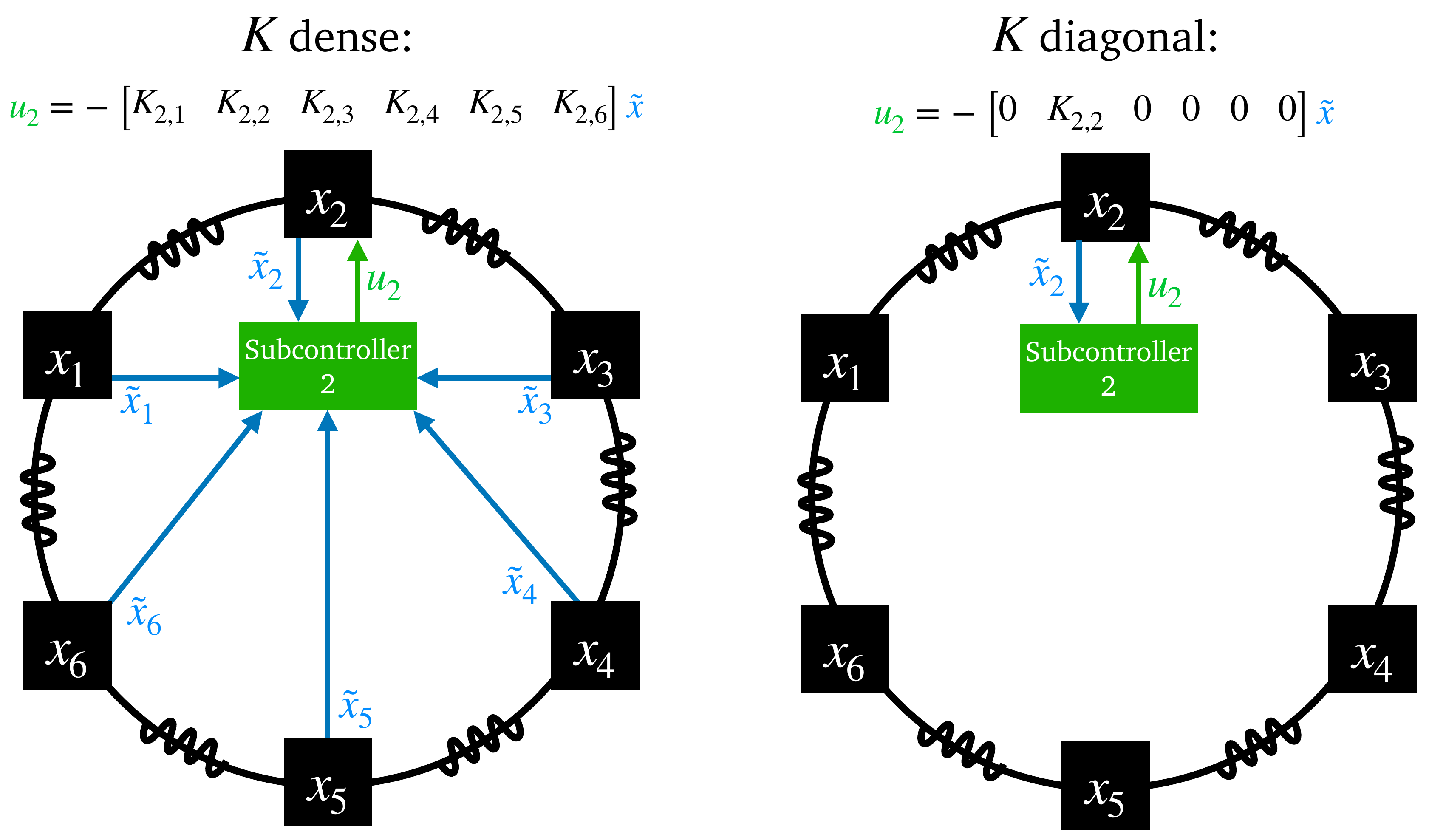}
    \caption{Mass-spring system in a circle with state $x$ and state estimate $\tilde{x}$. When the feedback gain $K$ is dense (left), computing the control $u$ at the second spatial location requires communication with all other locations. When $K$ is diagonal, or \textit{completely decentralized} (right) no communication is required.}
    \label{fig:commsavings}
\end{figure}

In this paper, we seek \textit{complete decentralization} of LQR and KF gains. This intuitively means the gains can be implemented completely locally, that measurements and state estimates need not be sent to other spatial locations. In finite dimensions, this means the gain matrices are comprised of only diagonal components (see, for example, Fig. \ref{fig:commsavings}). When both the KF and LQR gains are completely decentralized, we say the \textit{LQG controller is completely decentralized}. 

\section{Linear Quadratic Gaussian Control for Wave Equation} \label{sec:problemPDE}

Consider the following wave dynamics:
\begin{align}
\label{eq:wave}\partial_t^2p(x,t)&=c^2\partial_x^2p(x,t)+u(x,t)+d(x,t), &\\
    y(x,t)&=p(x,t)+m(x,t)\notag
 \end{align}
$\forall x\in\mathbb{T}=\left[-\frac{h}{2},\frac{h}{2}\right], \forall t\ge0$  where $p,u,y,d,m$ are the state, control, measurement, model disturbance, and measurement noise, respectively. The spatial domain $\mathbb{T}$ enforces periodic boundary conditions on all signals, i.e. $p\left(-\tfrac{h}{2},t\right)=p\left(\tfrac{h}{2},t\right)$  $\forall t$.
We write (\ref{eq:wave}) in state space form\footnote{This representation is not unique, but is the proper form for this problem as others are not able to uniquely capture measurements of $p$. Derivatives in (\ref{eq:ssLQG}) are defined in the weak sense to allow for mild solutions.} as 
\begin{align}\label{eq:ssLQG}\partial_t
 \begin{bmatrix}
   p(x,t)\\\partial_tp(x,t)\end{bmatrix}&=\underbrace{\begin{bmatrix}
      0&I\\c^2\partial_x^2&0
   \end{bmatrix}}_{\cA}\begin{bmatrix}
      p(x,t)\\\partial_tp(x,t)
 \end{bmatrix}+\underbrace{\begin{bmatrix}
   0\\I
 \end{bmatrix}}_{\cB} (u
+d)\\
y(x,t)&=\underbrace{\begin{bmatrix}
I &0
\end{bmatrix}}_{\cC}\begin{bmatrix}
    p(x,t)\\\partial_tp(x,t)
\end{bmatrix}+m.\notag\end{align}
$\cA$ generates a $C_0$-semigroup on $\chi_\alpha=\cH_\alpha(\mathbb{T})\oplus L_2(\T)$ for any $\alpha>0$ \cite{jensen2020localization, curtain2020introduction}, thus for each fixed $t$ we define \[P(\cdot,t)=\begin{bmatrix}
   p(\cdot,t)\\\partial_tp(\cdot,t)
\end{bmatrix}\in\chi_\alpha,\hspace{0.2in}u(\cdot,t)\in L_2(\T).\] 

 The objective of LQG is to design an optimal output feedback law, i.e. find $\cK$ and $\calL$ such that the state estimate $\tilde{P}$, defined by the dynamics 
\begin{align}\label{eq:LQGsolutionPDE}
    \partial_t{\tilde{P}}(x,t)&=(\cA-\calL\cC)\tilde{P}(x,t)+\cB u(x,t)+\calL y(x,t)\\
    u(x,t)&=-\cK\tilde{P}(x,t)\notag,
\end{align}
minimizes the quadratic cost function
\begin{equation}\label{eq:PDELQGcost}
    \cJ_{LQG}=\int_0^\infty\langle P, \cQ P\rangle_{\chi_\alpha}+\langle u, \cR u\rangle_{L_2(\T)}\,dt
\end{equation}
for positive definite operators $\cQ\in\cL(\chi_\alpha)$, $\cR\in\cL(L_2(\T))$. We choose $\cQ$ and $\cR$ to be of the form \begin{equation}R=\frac{1}{r^2} I ,\hspace{0.2in}Q=\begin{bmatrix}
   \frac{1}{q_1^2} I &0\\0&\frac{1}{q_2^2} I
   \end{bmatrix}, \hspace{0.2in}q_1,q_2,r>0.\end{equation}  

By the separation principle \cite{robust_text}, the optimal control policy is independent of the optimal estimator. That is, we can solve the LQR problem for the optimal $\cK$, and solve the Kalman filtering problem for $\calL$ and combine them in the LQG solution to retain optimality. We will solve the spatially discretized (and therefore finite dimensional) analogs of the above infinite dimensional LQR and KF problems, and analyze the locality of the gain matrices as a function of four nondimensional parameters. 

\section{LQG for the Discrete Wave Equation: A first attempt}\label{sec:LQGL2}
A discretization of the dynamics (\ref{eq:ssLQG}) in space represents an approximation to the full PDE dynamics, and also describes systems that are inherently discrete \cite{lindberg2024h,jovanovic2005ill}. The discretized analogs of the state, control, disturbance, and noise variables are given by $\bm{p}(t),\bm{u}(t),\bm{d}(t),\bm{m}(t)\in\R^n$. We approximate the second spatial derivative with periodic boundary conditions with a second-order Taylor approximation, i.e. $\partial_x^2\approx\frac{1}{(\Delta x)^2}\mathbb{D}^2$ 
where $\Delta x=\frac{h}{n}$ is the distance between discretization points and $\mathbb{D}^2$ is the $n\times n$ circulant matrix with first row $\begin{bmatrix}
    -2&1&0&\ldots&0&1
\end{bmatrix}$ so that  \[\lim_{\Delta x\to0}\frac{1}{\Delta x^2}\mathbb{D}^2=\lim_{n\to\infty}\frac{1}{\Delta x^2}\mathbb{D}^2=\partial_x^2.\]
Then the discretized analog of the dynamics (\ref{eq:ssLQG}) 
are 
   \begin{align}\label{eq:ssLQGdiscrete}
    \frac{d}{dt}{\bm{P}}(t)&=\underbrace{\begin{bmatrix}
         0&I\\\frac{c^2}{\Delta x^2}\mathbb{D}^2&0
      \end{bmatrix}}_{A}\bm{P}(t)+\underbrace{\begin{bmatrix}
      0\\I
    \end{bmatrix}}_{B} (\bm{d}(t)+\bm{u}(t))\\
    \bm{y}(t)&=\underbrace{\begin{bmatrix}
      I&0
    \end{bmatrix}}_C\bm{P}(t)+\bm{m}(t)\notag
   \end{align}
   where $\bm{P}(t)=\begin{bmatrix}
      \bm{p}(t)&\tfrac{d}{dt}{\bm{p}}(t)\end{bmatrix}^T$.

\begin{figure}
    \centering
\includegraphics[width=\linewidth]{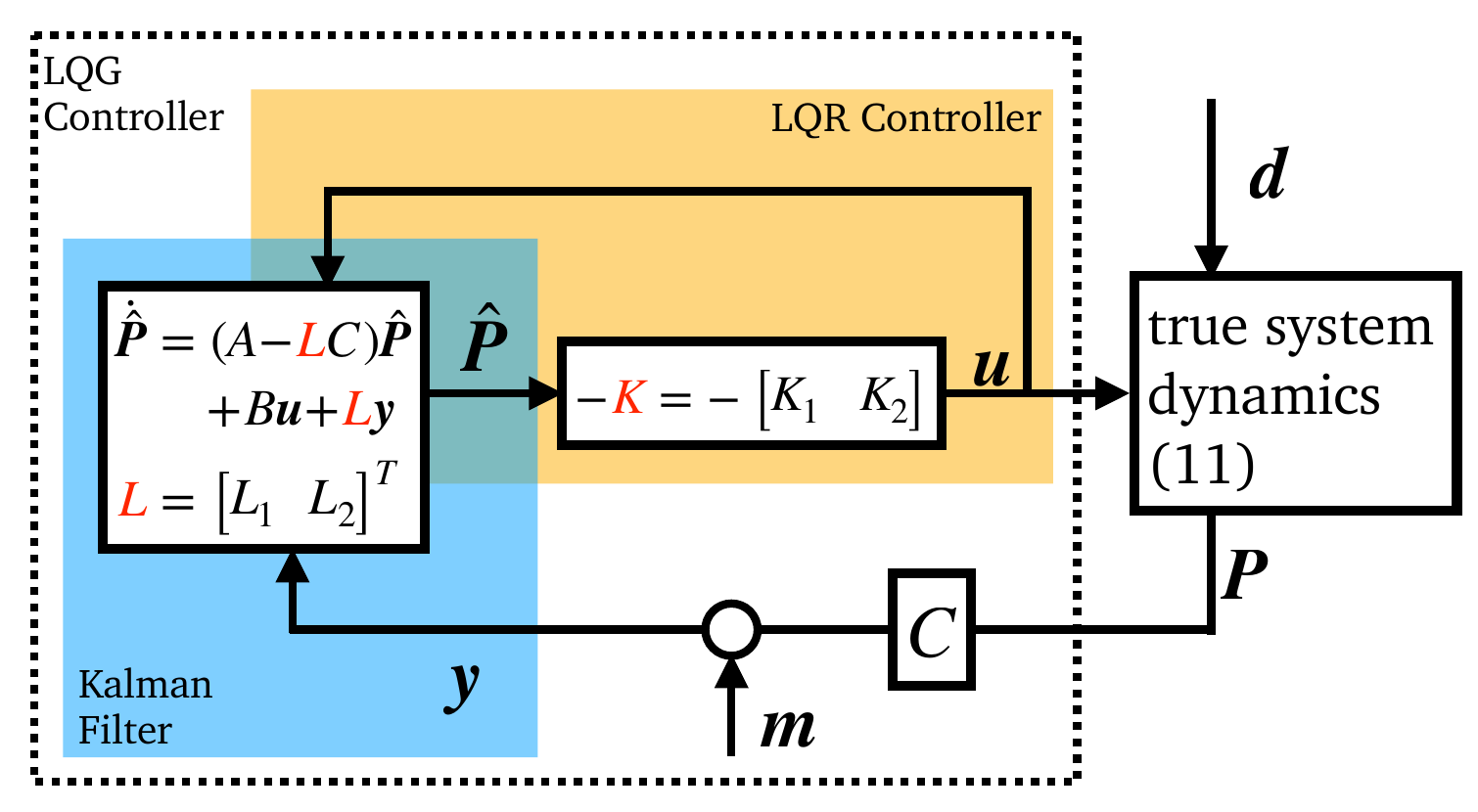}    \caption{Diagram of LQG controller. By the separation principle, LQG can be split into a Kalman filter (in blue) with gain $L$, and an LQR controller (in yellow) with gain $K$. In order to \textit{completely decentralize} $K$ and $L$ (in red), $K_1,K_2,L_1,L_2$ must all be diagonal matrices.}
    \label{fig:LQGdiagram}
\end{figure}

   We use finite dimensional linear systems theory to solve the LQG problem with dynamics (\ref{eq:ssLQGdiscrete}), visualized in Fig \ref{fig:LQGdiagram}.
   Lemmas 1-2 provide analytic expressions for the optimal LQR gain $K$ and KF gain $L$, and demonstrate that complete decentralization of the LQG problem is never optimal in this formulation. 

   \subsection{Linear Quadratic Regulator}
   We consider \eqref{eq:ssLQGdiscrete} without model disturbance and with full and perfect access to the state. The LQR problem is 
   \begin{subequations} \begin{align}
        &\min~ \int_0^\infty \bm{P}^T\begin{bmatrix}\frac{1}{q_1^2} I&0\\0&
            \frac{1}{q_2^2} I\end{bmatrix}\bm{P}+\frac{1}{r^2}\bm{u}^T\bm{u} \,dt \label{eq:costLQRL2}\\
   \label{eq:ssdiscLQR}
   & ~{\rm s.t.}~ ~\frac{d}{dt}{
    \bm{P}}(t)=\begin{bmatrix}
         0&I\\\frac{c^2}{\Delta x^2}\mathbb{D}^2&0
      \end{bmatrix}\bm{P}(t)+\begin{bmatrix}
      0\\I
    \end{bmatrix} \bm{u}(t).
   \end{align}\end{subequations}

      For ease of analysis, we nondimensionalize our system. The dimensions are time ($t$), state ($s$), and length $(\ell)$. We have parameters $r,q_1,q_2,\Delta x,c$ with dimensions
\[[r]=\frac{s}{t^2},\,[q_1]=s,\,[q_2]=\frac{s}{t}\,,\,\left[\Delta x\right]=\ell,\,[c]=\frac{\ell}{t}.\]
Define the nondimensional space, time, state, and control as 
\begin{equation}\label{eq:LQRnondim}\chi=\frac{x}{\Delta x}, \tau=\frac{ct}{\Delta x}, \bm{\phi}(\tau)=\frac{c^2}{r\Delta x^2}\bm{p}(\tau), \bm{\omega}=\frac{\bm{u}(\tau)}{r}.\end{equation} {We use $\Delta x$ instead of the interval length $h$ to avoid explicit dependence on  $n$ in the dynamics.}  
Then the nondimensional analog of \eqref{eq:costLQRL2}-\eqref{eq:ssdiscLQR} is 
\begin{subequations} \begin{align}
 & \min~  \int_0^\infty \bm{\Phi}^T\begin{bmatrix}I&0\\0&\Pi_2I\end{bmatrix}\bm{\Phi}+\frac{1}{\Pi_3^2}\bm{\omega}^T\bm{\omega} \,d\tau \label{eq:nondimLQRcostL2} \\
 &~~ {\rm s.t.} ~\frac{d}{d\tau}{\bm{\Phi}}(\tau)
   =\begin{bmatrix}
        0&I\\\mathbb{D}^2&0
     \end{bmatrix}\bm{\Phi}(\tau)+\begin{bmatrix}
     0\\I
   \end{bmatrix}\bm{\omega}(\tau) \label{eq:ssnondimLQR}
  \end{align} 
  \end{subequations}
where $\bm{\Phi}(\tau)=\begin{bmatrix}
    \bm{\phi}(\tau)&\tfrac{d}{d\tau}{\bm{\phi}}(\tau)
\end{bmatrix}^T$ and $\Pi_2,\Pi_3$ are the nondimensional parameters \begin{equation}\label{eq:pi2pi3}\Pi_2=\frac{c^2q_1^2}{q_2^2\Delta x^2},\Pi_3=\frac{\Delta x^2r}{c^2q_1}.\end{equation}

\begin{lemma}\label{thm:LQRsolutionL2}
    The optimal feedback gain $K$ for the LQR problem  (\ref{eq:ssnondimLQR}, \ref{eq:nondimLQRcostL2}) can be represented in the nondimensional spatial frequency domain (with $\kappa$ denoting the nondimensional spatial frequency) as 
    \begin{equation}\label{eq:finalKL2}
        \hat{K}(\kappa)=\begin{bmatrix}
            \hat{K}_0(\kappa)&\sqrt{2\hat{K}_0(\kappa)+{\Pi_2}{\Pi_3^2}}
         \end{bmatrix}
    \end{equation}
    where \begin{equation}\hat{K}_0(\kappa)=\hat{D}_{\kappa\kappa}+\sqrt{\hat{D}^2_{\kappa\kappa}+\Pi_3^2}\end{equation} and $\hat{D}_{\kappa\kappa}=-4\sin^2\left(\frac{\pi \kappa}{n}\right)$. K can be recovered as $K=\begin{bmatrix}
        K_1&K_2\end{bmatrix}$ where $K_1$ and $K_2$ are the circulant matrices 
    \begin{align}
        K_1&=F^{-1}\mathrm{diag}\left(\hat{K_0}(\kappa)\right)F\\
        K_2&=F^{-1}\mathrm{diag}\left(\sqrt{2\hat{K}_0(\kappa)+  {\Pi_2}{\Pi_3^2}}\right)F.
    \end{align}
\end{lemma}
\begin{proof}
    See Appendix \ref{sec:thm1proof}, setting $\Pi_1=0$.
\end{proof}

\subsection{Estimation via Kalman Filter}
We now consider the wave dynamics in (\ref{eq:ssLQGdiscrete}) without control. The corresponding estimation problem\footnote{Here $\bm{m},\bm{d}$ are deterministic but unknown signals. This deterministic formulation of the Kalman filter is somewhat atypical, but allows us to avoid making assumptions on the distribution of noise in the system. Minimizing (\ref{eq:costdiscKFL2}) is equivalent to minimizing the asymptotic estimation error, $\lim_{t\to\infty}\mathbb{E}[\|\bm{P}(t)-\hat{\bm{P}}(t)\|^2]$,
when $\bm{m}$ and $\bm{d}$ are uncorrelated zero-mean
Gaussian white-noise stochastic processes with covariance matrices $\delta(t-\bar{t})\sigma^2_m I$ and $\delta(t-\bar{t})\sigma_d^2I$, respectively, where $\delta$ refers to a Dirac delta. See \cite[Section~23.2-3]{hespanha2018linear} and \cite{willems2004deterministic} for more details.} is 
\begin{subequations} \begin{align}
        \min~ &\int_0^\infty \frac{1}{\sigma_m^2}\bm{m}^T\bm{m} +\frac{1}{\sigma_d^2}\bm{d}^T\bm{d} \,dt \label{eq:costdiscKFL2}\\
   \label{eq:ssestdiscrete}
    ~{\rm s.t.}~& ~
    \frac{d}{dt}{\bm{P}}(t)=\begin{bmatrix}
         0&I\\\frac{c^2}{\Delta x^2}\mathbb{D}^2&0
      \end{bmatrix}\bm{P}(t)+\begin{bmatrix}
      0\\I
    \end{bmatrix} \bm{d}(t)\\
    &\hspace{0.26in}\bm{y}(t)=\begin{bmatrix}
      I&0
    \end{bmatrix}\bm{P}+\bm{m}(t)\notag.
   \end{align}\end{subequations}
for constants $\sigma_m,\sigma_d>0$ where the ratio of $\sigma_m$ and $\sigma_d$ quantifies the ratio of model accuracy to sensor quality.

We again nondimensionalize the system for analysis. 
 We have parameters $\sigma_m,\sigma_d,\Delta x,c$ with dimensions
\[[\sigma_d]=\frac{s}{t^2},\,[\sigma_m]=s,\,\left[\Delta x\right]=\ell,\,[c]=\frac{\ell}{t}.\]
Define the nondimensional space, time, state, control, disturbance, noise, and measurement variables as 
\begin{align}\label{eq:KFnondim}
   \chi=\frac{x}{\Delta x},\tau=\frac{ct}{\Delta x},\bm{\psi}(\tau)=\frac{c^2}{\Delta x^2\sigma_d}\bm{p}(\tau),\\ \bm{\rho}=\frac{\bm{d}(\tau)}{\sigma_d}, \bm{\eta}=\frac{\bm{m}(\tau)}{\sigma_m},\bm{\gamma}(\tau)=\frac{1}{\sigma_m}\bm{y}(\tau).\notag
\end{align} 
{This nondimensionalization is different than (\ref{eq:LQRnondim}) since the parameters $\sigma_d$ and $\sigma_m$ may differ from LQR weights $r, q_1,q_2$.} The nondimensional problem is 
\begin{subequations} \begin{align}
        \min~ &\int_0^\infty \bm{\eta}^T\bm{\eta}+\bm{\rho}^T\bm{\rho} \,d\tau \label{eq:nondimKFcostL2}\\
   \label{eq:ssnondimKF}
    ~{\rm s.t.}~& ~\frac{d}{d\tau}{
  \bm{\Psi}}(\tau)=\begin{bmatrix}
        0&I\\\mathbb{D}^2&0
     \end{bmatrix}\bm{\Psi}(\tau)+\begin{bmatrix}
     0\\I
   \end{bmatrix} \bm{\rho}(\tau)\\
   &\hspace{0.28in}\bm{\gamma}(\tau)=\begin{bmatrix}
      \Pi_4 & 0
   \end{bmatrix}\bm{\Psi}(\tau)+\bm{\eta}(\tau)\notag.
   \end{align}\end{subequations}
  where $\bm{\Psi}(\tau)=\begin{bmatrix}
      \bm{\psi}(\tau) &\tfrac{d}{d\tau}{\bm{\psi}}(\tau)
 \end{bmatrix}^T$ and $\Pi_4$ is the nondimensional parameter $\Pi_4=\frac{\Delta x^2\sigma_d}{c^2\sigma_m}$.  

\begin{lemma}\label{thm:KFsolutionL2}
   The optimal Kalman filter gain $L$ for the estimation problem described (\ref{eq:ssnondimKF}, \ref{eq:nondimKFcostL2}) can be represented in the nondimensional spatial frequency domain as 
   \begin{equation}\label{eq:finalLL2}
       \hat{L}(\kappa)=\begin{bmatrix}
         \sqrt{\frac{2}{\Pi_4}\hat{L}_0(\kappa)}& \hat{L}_0(\kappa)
        \end{bmatrix}^T
   \end{equation}
   where \begin{equation}\hat{L}_0(\kappa)=\frac{\hat{D}_{\kappa\kappa}}{\Pi_4}+\sqrt{\frac{\hat{D}^2_{\kappa\kappa}}{\Pi_4^2}+1}\end{equation} and $\hat{D}_{\kappa\kappa}=-4\sin^2\left(\frac{\pi \kappa}{n}\right)$. L can be recovered as $L=\begin{bmatrix}
        L_1&L_2\end{bmatrix}^T$ where $L_1$ and $L_2$ are  the circulant matrices
    \begin{align}
        L_1&=F^{-1}\mathrm{diag}\left(\sqrt{\frac{2}{\Pi_4}\hat{L}_0(\kappa)}\right)F\\
        L_2&=F^{-1}\mathrm{diag}\left(\hat{L_0}(\kappa)\right)F.
    \end{align}
\end{lemma}
\begin{proof}
    See Appendix \ref{sec:thm2proof}, setting $\Pi_1=0$.
\end{proof}
\subsection{Inability to decentralize}\label{sec:issue}
Complete decentralization requires that the gain matrices be constant over the nondimensional spatial frequency $\kappa$. We observe that no choice of $\Pi_3$ or $\Pi_4$ makes $\hat{K}_0(\kappa)$ constant in Lemma \ref{thm:LQRsolutionL2} or $\hat{L}_0(\kappa)$ constant in Lemma \ref{thm:KFsolutionL2}, therefore complete decentralization is not optimal for any choice of parameters. However, the cost functions \eqref{eq:nondimLQRcostL2}, \eqref{eq:nondimKFcostL2} use $\langle \cdot,\cdot \rangle_{\ell_2}$, which is the discrete analog of $\langle \cdot,\cdot \rangle_{L_2}$, but the PDE requires that $P(x,t)\in\chi_\alpha$, reflected in the inner product in (\ref{eq:PDELQGcost}). In the following section, we align the discrete cost functions with the PDE cost functions, and show that this allows for complete decentralization of the optimal LQG controller.

\section{PDE-Informed LQG }\label{sec:smartLQG}
We return to the full PDE LQG problem to motivate new cost functions for the spatially discretized LQG controller. 

\subsection{Linear Quadratic Regulator}
The LQR problem for the full PDE \eqref{eq:ssLQG} finds the spatio-temporal signal $u$ that minimizes 
   \begin{align}\label{eq:firstLQRcost}
      \cJ_{LQR}&=\int_0^\infty \langle P,\cQ P\rangle_{\chi_\alpha}+\langle u ,\cR u\rangle_{L_2} \, dt,
   \end{align} where $P\in\chi_\alpha,u\in L_2$, subject to dynamics
\begin{align}\label{eq:ssLQR}\partial_t
 P(x,t)&=\cA P(x,t)+\cB u(x,t).
\end{align}
    On $L_2(\T)$, the adjoint of the operator $\partial_x$ is $-\partial_x$ so that
   \begin{align}
      \cJ_{LQR}
      \label{eq:PDEcostLQR}&=\int_0^\infty \Big\langle p,\frac{1}{q_1^2} p\Big\rangle_{L_2} -\alpha^2\Big\langle p,\frac{1}{q_1^2} \partial_x^2p\Big\rangle_{L_2}\notag\\
      & \hspace{0.5in}+ \Big\langle \partial_tp,\frac{1}{q_2^2} \partial_tp\Big\rangle_{L_2}+\Big\langle u ,\frac{1}{r^2}u\Big\rangle_{L_2} dt.
      \end{align}

   The corresponding discretized cost function is
      \begin{align}
         \int_0^\infty \bm{P}^T\begin{bmatrix}\frac{1}{q_1^2}(I- \frac{\alpha^2}{\Delta x^2}\mathbb{D}^2)&0\\0&
            \frac{1}{q_2^2}\end{bmatrix}\bm{P}+\frac{1}{r^2}\bm{u}^T\bm{u} \,dt.
      \end{align} 
      We nondimensionalize the same as in (\ref{eq:LQRnondim}) but now we have the additional parameter $\alpha$ with dimension $[\alpha]=\ell$. The dynamics remain the same as (\ref{eq:ssnondimLQR}) but with cost function 
  \begin{align} \label{eq:nondimLQRcost}
   J_{LQR}&=\int_0^\infty \bm{\Phi}^T\begin{bmatrix}I- \Pi_1\mathbb{D}^2&0\\0&\Pi_2I\end{bmatrix}\bm{\Phi}+\frac{1}{\Pi_3^2}\bm{\omega}^T\bm{\omega} \,d\tau
\end{align} 
where $\Pi_1=\frac{\alpha^2}{\Delta x^2}$ and $\Pi_2,\Pi_3$ remain as in \eqref{eq:pi2pi3}. 

Recall that $-\bm{\phi}^T\mathbb{D}^2\bm{\phi}$ is the discrete analog of the inner product $-\langle p,\partial_x^2p\rangle=\langle \partial_x p,\partial_xp\rangle$,
thus $\Pi_1$ functions as a weight on the potential energy of the system. Likewise, $\Pi_2$ weights kinetic energy, and $1/\Pi_3$ weights control effort.

\begin{lemma}\label{thm:LQRsolution}
    The optimal feedback gain $K$ for the LQR problem (\ref{eq:ssnondimLQR}, \ref{eq:nondimLQRcost}) can be represented in the nondimensional spatial frequency domain as 
    \begin{equation}\label{eq:finalK}
        \hat{K}(\kappa)=\begin{bmatrix}
            \hat{K}_0(\kappa)&\sqrt{2\hat{K}_0(\kappa)+{\Pi_2}{\Pi_3^2}}
         \end{bmatrix}
    \end{equation}
    where \begin{equation}\hat{K}_0(\kappa)=\hat{D}_{\kappa\kappa}+\sqrt{\hat{D}^2_{\kappa\kappa}+\Pi_3^2\left(1-\Pi_1\hat{D}_{\kappa\kappa}\right)}\end{equation} and $\hat{D}_{\kappa\kappa}=-4\sin^2\left(\frac{\pi \kappa}{n}\right)$. K can be recovered as $K=\begin{bmatrix}
        K_1&K_2\end{bmatrix}$ where $K_1$ and $K_2$ are the circulant matrices 
    \begin{align}
        K_1&=F^{-1}\mathrm{diag}\left(\hat{K_0}(\kappa)\right)F\\
        K_2&=F^{-1}\mathrm{diag}\left(\sqrt{2\hat{K}_0(\kappa)+{\Pi_2}{\Pi_3^2}}\right)F.
    \end{align}
\end{lemma}
\proof{
    See Appendix \ref{sec:thm1proof}.
\qed}

When $\Pi_1=\frac{2}{\Pi_3}$, $ \hat{K}_0(\kappa)=\Pi_3\,\, \forall \kappa$, so that $K_1$ and $K_2$ are diagonal matrices, i.e. $K$ is \textit{completely decentralized}. 
The condition $\Pi_1=\frac{2}{\Pi_3}$ is a relationship between the weights on potential energy and control effort, and is independent of the kinetic energy weight. Complete decentralization is impossible when $\Pi_1=0$, which is consistent with our result in Section \ref{sec:issue}. Having nonzero $\Pi_1$ allows us to consider potential energy in our system and is equivalent to incorporating a `Sobolev-like' norm on the state $\phi$, as opposed to an $L_2$ norm that ignores the spatial derivative. 
 Including potential energy in our cost not only makes physical sense but also allows for the desired locality.

\subsection{Estimation via Kalman Filter}
 The wave dynamics without control for the PDE in (\ref{eq:ssLQG}) are
\begin{align}\label{eq:ssKFPDE}\partial_t
 P(x,t)&=\begin{bmatrix}
      0&I\\c^2\partial_x^2&0
   \end{bmatrix}P(x,t)+\begin{bmatrix}
   0\\I
 \end{bmatrix} d(x,t)\\
 y(x,t)&=\begin{bmatrix}
   I&0
 \end{bmatrix}P(x,t)+ m(x,t)\notag
\end{align}   
$\forall x\in\T ,\forall t\ge0$. The dynamics suggest we assume $d(x,t)\in L_2(\T)$ and $ m(x,t)\in \cH_\alpha(\mathbb{T})$.

   The Kalman filter chooses the estimate of $p(x,t)$ that minimizes 
   \begin{align}
      \label{eq:costKFPDE}
      \cJ_{KF}&=\int_0^\infty \Big\langle m,\frac{1}{\sigma_m^2}m\Big\rangle_{\cH_\alpha}+\Big\langle  d,\frac{1}{\sigma_d^2}d\Big\rangle_{L_2} dt
   \end{align} subject to dynamics (\ref{eq:ssKFPDE}) for constants 
   $\sigma_m,\sigma_d>0$. 
   The discretized dynamics are still given by (\ref{eq:ssestdiscrete})
  with the discrete analog of the cost function (\ref{eq:costKFPDE}) given by
      \begin{align}\label{eq:costdiscKF}
         \int_0^\infty \bm{m}^T\left(\frac{1}{\sigma_m^2}I-\frac{\alpha^2}{\sigma_m^2\Delta x^2}\mathbb{D}^2\right)\bm{m} +\frac{1}{\sigma_d^2}\bm{d}^T\bm{d} \,dt.
      \end{align} 
 We nondimensionalize the same as (\ref{eq:KFnondim}) so that dynamics remain the same as (\ref{eq:ssnondimKF}) but now the cost is 
  \begin{align}\label{eq:nondimKFcost}
   J_{KF}&=\int_0^\infty \bm{\eta}^T\left(I- \Pi_1\mathbb{D}^2\right)\bm{\eta}+\bm{\rho}^T\bm{\rho} \,d\tau
\end{align} 
where again $\Pi_1=\frac{\alpha^2}{\Delta x^2}$. In the equivalent stochastic formulation of the Kalman filter, according to \eqref{eq:nondimKFcost}, $\bm{\eta}$ and $\bm{\rho}$ are assumed to be Gaussian white noise with inverse covariance matrices $\cR_{\bm{\eta}}^{-1}=(I-\Pi_1\mathbb{D}^2)$ and $\cR_{\bm{\rho}}^{-1}=I$, respectively, i.e. 
\begin{equation}\E[\bm{\rho}(t)^T\bm{\rho}(\tau)]=\delta(t-\tau)\cdot I, \,\,\,\E[\bm{\eta}(t)^T\bm{\eta}(\tau)]=\delta(t-\tau)\cdot \cR_{\bm{\eta}}.\notag\end{equation} 
 We focus on $\bm{\eta}$. With $\Pi_1=0$, measurement noise is uncorrelated in space, but for $\Pi_1>0$ the noise is completely correlated with strength of correlation decaying with distance, which is likely more realistic. 
\begin{lemma}\label{thm:KFsolution}
   The optimal Kalman filter gain $L$ for the estimation problem (\ref{eq:ssnondimKF}, \ref{eq:nondimKFcost}) can be represented in the nondimensional spatial frequency domain as 
   \begin{equation}\label{eq:finalL}
       \hat{L}(\kappa)=\begin{bmatrix}
         \sqrt{\frac{2}{\Pi_4}\hat{L}_0(\kappa)}& \hat{L}_0(\kappa)
        \end{bmatrix}^T
   \end{equation}
   where \begin{equation}
   \hat{L}_0(\kappa)=\frac{\hat{D}_{\kappa\kappa}}{\Pi_4}+\sqrt{\frac{\hat{D}^2_{\kappa\kappa}}{\Pi_4^2}-\Pi_1\hat{D}_{\kappa\kappa}+1}\end{equation}  
   and $\hat{D}_{\kappa\kappa}=-4\sin^2\left(\frac{\pi \kappa}{n}\right)$. L can be recovered as $L=\begin{bmatrix}
        L_1&L_2\end{bmatrix}^T$ where $L_1$ and $L_2$ are the circulant matrices
    \begin{align}
        L_1&=F^{-1}\mathrm{diag}\left(\sqrt{\frac{2}{\Pi_4}\hat{L}_0(\kappa)}\right)F\\
        L_2&=F^{-1}\mathrm{diag}\left(\hat{L_0}(\kappa)\right)F.
    \end{align}
\end{lemma}
\begin{proof}
    See Appendix \ref{sec:thm2proof}. 
\end{proof}

When $\Pi_1=\frac{2}{\Pi_4}$, $\hat{L}_0(\kappa)=1\,\,  \forall \kappa $ so that the Kalman filter gain $L$ is  \textit{completely decentralized}. The condition $\Pi_1=\frac{2}{\Pi_4}$ is impossible when $\Pi_1=0$, thus accounting for spatial correlations in measurement noise in the Kalman filter cost may be more realistic and also allow for the desired localization results. 

Combining the results from Lemma \ref{thm:LQRsolution} and Lemma \ref{thm:KFsolution}, and choosing $\sigma_d=r$ so that $\bm{\Psi}=\bm{\Phi}$, we have the optimal nondimensional discrete LQG controller 
\begin{align}
    \frac{d}{d\tau}{\tilde{\bm{\Phi}}}(\tau)&=(A-LC)\tilde{\bm{\Phi}}(\tau)+B\bm{\omega}(\tau)+L\bm{\gamma}(\tau)\\
    \bm{\omega}&=-K\tilde{\bm{\Phi}}(\tau)\notag.
\end{align}
Computing $\tilde{\bm{\Phi}}(\tau)$ requires spatial communication via the dynamics $A-LC$, but $A$ only requires communication with two nearest neighbors, thus localizing $K$ and $L$ can indeed reduce communication costs.

\section{LQG Locality and Performance Analysis} \label{sec:discussion} 

We have shown that reformulating the LQR and Kalman filtering problem with `Sobolev-like' cost functions allows for complete decentralization of both gain matrices. We now give a condition for complete decentralization of the LQG problem, and characterize the relationship between controller locality and performance. 

\begin{theorem}
    The optimal LQG controller for the nondimensional discrete wave equation (described in Section~\ref{sec:smartLQG}) is completely decentralized when the LQG weights are chosen so that $q_1=\sigma_m$ and $r=\sigma_d$, and $\Pi_1=\frac{2}{\Pi_4}$.
\end{theorem}
\begin{proof}
First note when $q_1=\sigma_m$ and $r=\sigma_d$, $\Pi_3=\Pi_4$. Let $\Pi_1=\frac{2}{\Pi_4}=\frac{2}{\Pi_3}$. Since $\Pi_1=\frac{1}{\Pi_3}$, by Lemma \ref{thm:LQRsolution} we have $\hat{K}_0=\Pi_3$ for all $\kappa$. Then 
\begin{equation}
    \hat{K}(\kappa)=\begin{bmatrix}
        \Pi_3 & \sqrt{2\Pi_3+{\Pi_2}{\Pi_3^2}}
    \end{bmatrix}
\end{equation} so that 
    $K=\begin{bmatrix}
        \mathrm{diag}(\Pi_3) & \mathrm{diag}\left(\sqrt{2\Pi_3+{\Pi_2}{\Pi_3^2}}\right)
    \end{bmatrix}.$
Similarly, since $\Pi_1=\frac{2}{\Pi_4}$, by Lemma \ref{thm:KFsolution} we have $\hat{L}_0=1$ for all $\kappa$. Then  
\begin{equation}
    \hat{L}(\kappa)=\begin{bmatrix}
        \sqrt{\frac{2}{\Pi_4}}&1
    \end{bmatrix}^T
\end{equation} so that 
    $L=\begin{bmatrix}
        \mathrm{diag}\left(\sqrt{\frac{2}{\Pi_4}}\right) &\mathrm{diag}(1)
    \end{bmatrix}^T.$

\end{proof} 
In other words, appropriately choosing LQG weights and system parameters allows for significant communication savings, since measurement corrections and control computations can be performed completely locally.

    Assume we have chosen $\Pi_3=\Pi_4$ so that we are able to completely decentralize the LQG controller. Then the condition for complete locality is given by 
\begin{equation}\label{eq:localparams}\frac{\alpha^2\sigma_d}{c^2\sigma_m}=2.\end{equation} If the system is fixed, then $\alpha$ functions
as a tuning parameter to decentralize the system. If considering locality in system design, in addition to $\alpha$, the wavespeed $c$ can be altered (e.g. via choice of material), and $\sigma_d$ and $\sigma_m$ can be tuned as a function of model and sensor quality. The condition (\ref{eq:localparams}) is independent of $\Delta x$ which means decentralization does not depend on interval length $h$, nor discretization resolution $n$. 

\begin{figure}
\centering
\subfigure[]{
    \includegraphics[width=0.8\linewidth]{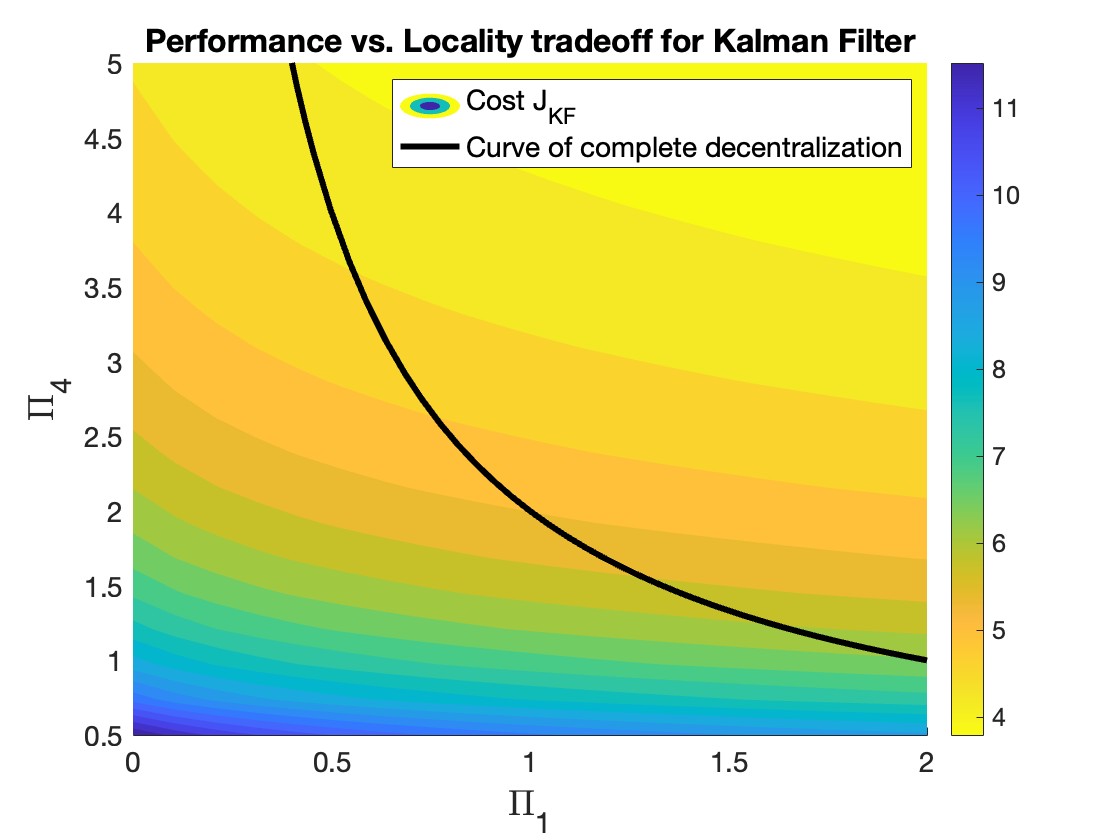}\label{fig:KFperf}
}
\subfigure[]{
    \includegraphics[width=0.8\linewidth]{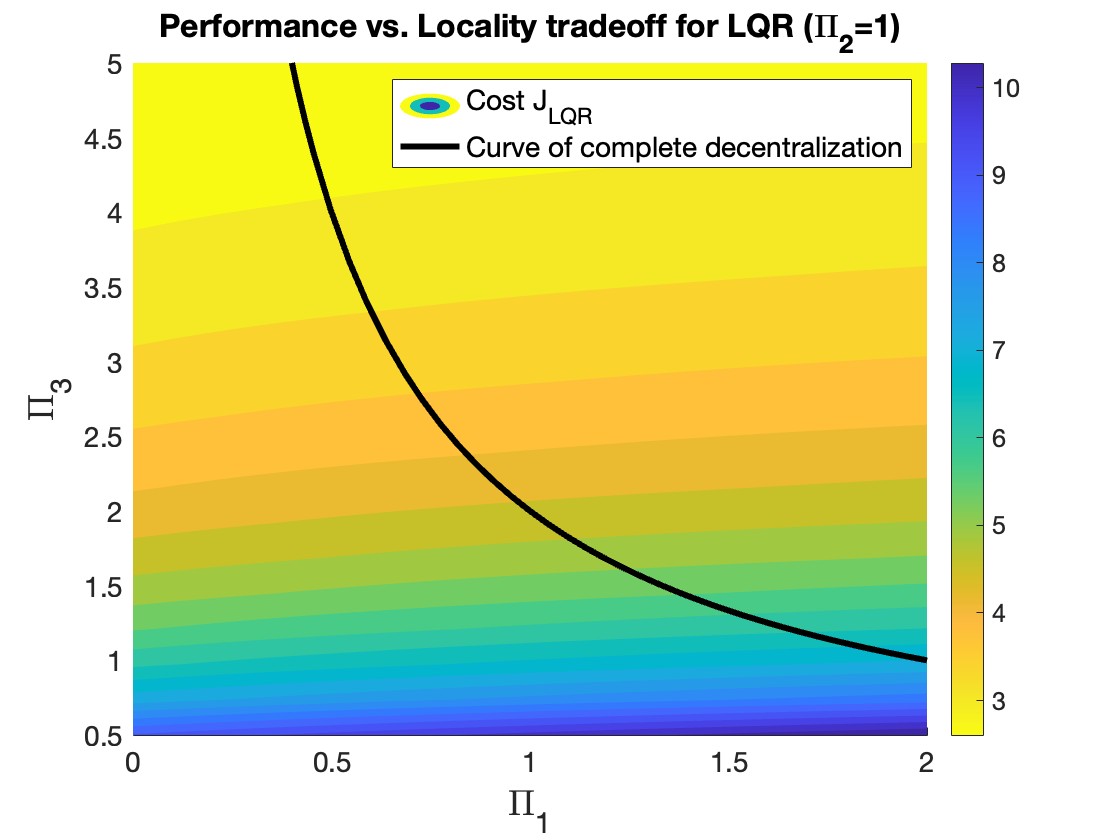}\label{fig:LQRperf}
}
 
   \caption{Tradeoff between locality and performance for the Kalman filter (Fig. \ref{fig:KFperf}) and LQR controller (Fig. \ref{fig:LQRperf}) with $n=30$. The best performance is in the yellow regions and complete decentralization lies on the black curve. We see that performance does not degrade with locality.}
   \label{fig:indivperf}
\end{figure}

\begin{figure}[]\centering
   \includegraphics[width =0.8\linewidth]{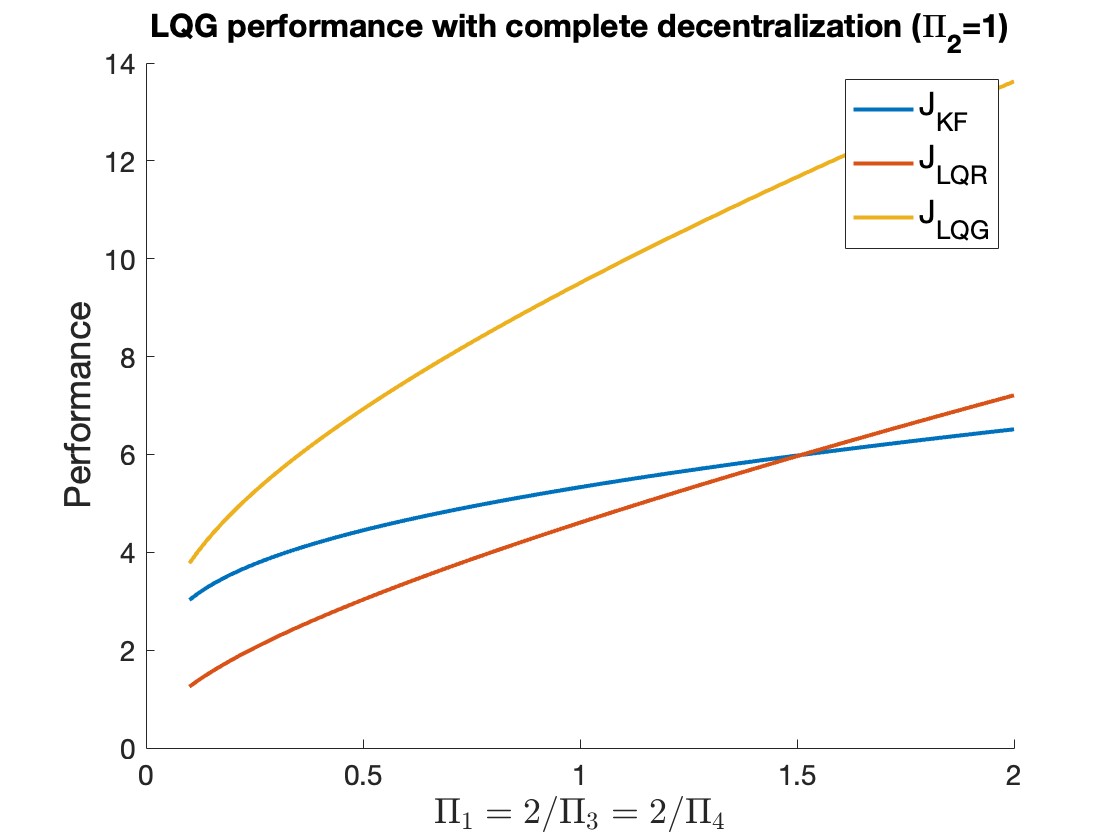}
   \caption{Performance along the black curves in Fig. \ref{fig:indivperf} and the corresponding total LQG cost. Pick $\Pi_1$ small to minimize cost while maintaining complete decentralization.}
   \label{fig:localperf}
\end{figure}

{It is interesting that \cite{arbelaiz2022information} found the same condition \eqref{eq:localparams} for complete decentralization of the Kalman filter for wave equation PDE dynamics on $\R$, utilizing spatial decay rates on an infinite domain.\footnote{A similar `matching condition' is derived in \cite{arbelaiz2024optimal} for estimating general PDE dynamics mapping $L_2(\mathbb{R})$ to itself. However, this result does not include the wave equation, which is not well posed over $L_2$.} We derive this condition for discretized dynamics on a compact domain where decay rates are not well defined, and we extend to the full LQG problem. 
Our finite-dimensional approach also allows us to quantify controller performance as a function of controller localization. LQR, KF, and LQG cost can be computed exactly via trace of appropriate matrices (see \cite{robust_text})}. In Fig. \ref{fig:indivperf} we explore how the Kalman filter and LQR systems perform as we vary the parameters $\Pi_1$ and $\Pi_4,\Pi_3$. {LQR and KF costs do not rise along the curve of decentralization, meaning there is not an inherent tradeoff between policy locality and performance. Rather, as shown in Fig. \ref{fig:localperf}, system parameters can be chosen to minimize the overall LQG cost while maintaining complete decentralization.}

\section{Conclusion}\label{sec:conclusion}
We have demonstrated that drawing inspiration from PDEs in discrete systems can have very beneficial results. We derived a condition for complete decentralization of an LQG controller for the discrete wave equation on the unit circle, and provided a physical interpretation of this condition. We also showed that controller localization does not degrade performance.

This result is a first step toward analysis of more complicated systems, such as incorporating damping or nonlinearities, or considering different boundary conditions, possibly employing an ``embedding approach" \cite{epperlein2016spatially}. We would also like to explore whether decentralization of other linear PDEs over Sobolev spaces depends on a ``Sobolev-like" discretization. Finally, we are interested in analyzing the full PDE dynamics with a discrete number of sensors and actuators, in line with the work of \cite{morris}, {and analyzing performance of our designed discrete controller for the infinite-dimensional system.} Characterizing robustness of structure or performance of the LQG controller analyzed here is also an open question. 

\bibliography{references}
\bibliographystyle{ieeetr}

\section{Appendix A}\label{sec:appA}
\subsection{Proof of Lemmas \ref{thm:LQRsolutionL2},\ref{thm:LQRsolution}}\label{sec:thm1proof}
\begin{proof}
   Plancherel's Theorem \eqref{eq:planch} motivates solving in Fourier space. The DFT of (\ref{eq:ssnondimLQR}, \ref{eq:nondimLQRcost}) gives $n$ decoupled $2\times2$ optimization problems:
   \begin{subequations} \begin{align}
 \min~  &\label{eq:nondimLQRcostfourier}
                \int_0^\infty \hat{\bm{\Phi}}_k^T\underbrace{\begin{bmatrix}1- \Pi_1\hat{D}_{kk}&0\\0&\Pi_2\end{bmatrix}}_{\hat{Q}_k}\hat{\bm{\Phi}}_k+\underbrace{\frac{1}{\Pi_3^2}}_{\hat{R}}\hat{\bm{\omega}}_k^T\hat{\bm{\omega}}_k \,d\tau \\
 ~~ {\rm s.t.} ~\label{eq:LQRfourierss}&\dot{\hat{\bm{\Phi}}}_k=\begin{bmatrix}
          0&1\\\hat{D}_{kk}&0
          \end{bmatrix}\hat{\bm{\Phi}}_k+
          \begin{bmatrix}0\\1\end{bmatrix}\hat{\bm{\omega}}_k:=\hat{A}_k\hat{\bm{\Phi}}
          _k+\hat{B}\hat{\bm{\omega}}_k
  \end{align} 
  \end{subequations}
             where $\hat{D}_{kk}=-4\sin^2\left(\frac{\pi k}{n}\right)$ are the eigenvalues of $\mathbb{D}^2$ and $\hat{\bm{\Phi}}_k=[\hat{\bm{\phi}}_k ~\dot{\hat{\bm{\phi}}}_k
            ]^T$. The solution to (\ref{eq:nondimLQRcostfourier}, 
             \ref{eq:LQRfourierss})  is given by
             $\hat{\bm{\omega}}_k=\hat{K}_k\hat{\bm{\Phi}}_k$ where \begin{align}
                \hat{K}_k=\hat{R}^{-1}\hat{B}^T\hat{P}_k&=\Pi_3^2\begin{bmatrix}
                    0&1
                 \end{bmatrix}\begin{bmatrix}
                 \hat{P}_1(k)&\hat{P}_0(k)\\\hat{P}_0(k)&\hat{P}_2(k)
                 \end{bmatrix}\label{eq:genK}
             \end{align}
             and $\hat{P}_k\succ 0$ is the solution to the Algebraic Riccati Equation 
             \begin{equation}\label{eq:LQRARE}
                \hat{P}_k\hat{A}_k+\hat{A}^T_k\hat{P}_k-\hat{P}_k\hat{B}\hat{R}^{-1}\hat{B}^T\hat{P}_k+\hat{Q}_k=0.
             \end{equation}
          Solving (\ref{eq:LQRARE}) gives 
          \begin{align}\hat{P}_0(k)&=\frac{1}{\Pi_3^2}\left(\hat{D}_{kk}+\sqrt{\hat{D}^2_{kk}+\Pi_3^2\left(1-\Pi_1\hat{D}_{kk}\right)}\right)\\
          \hat{P}_2(k)&=\frac{1}{\Pi_3}\sqrt{2\hat{P}_0(k)+\Pi_2}\end{align} so that (\ref{eq:genK}) can be written as (\ref{eq:finalK}).
\end{proof}
\subsection{Proof of Lemmas \ref{thm:KFsolutionL2},\ref{thm:KFsolution}}
\begin{proof}\label{sec:thm2proof}
   The DFT of (\ref{eq:ssnondimKF}, \ref{eq:nondimKFcost}) gives
    $n$ $2\times 2$ decoupled optimization problems:
    \begin{subequations} \begin{align}
 \min~  &\label{eq:nondimKFcostfourier}
               \int_0^\infty \hat{\bm{\eta}}^T_k{\hat{Q}_k}\hat{\bm{\eta}}_k+\hat{\bm{\rho}}_k^T\hat{\bm{\rho}} _k\,d\tau,~~ \hat{Q}_k := {1- \Pi_1\hat{D}_{kk}} \\
 ~~ {\rm s.t.} ~&\label{eq:KFfourierss}\dot{\hat{\bm{\Psi}}}_k=\begin{bmatrix}
         0&1\\\hat{D}_{kk}&0
         \end{bmatrix}\hat{\bm{\Psi}}_k+\begin{bmatrix}
         0\\1
         \end{bmatrix}\hat{\bm{\rho}}_k:=\hat{A}_k\bm{\Psi}_k+\hat{B}\hat{\bm{\rho}}_k\notag\\
         &\hspace{0.05in}\hat{\bm{\gamma}}_k=\begin{bmatrix}
         \Pi_4&0\end{bmatrix}\hat{\bm{\Psi}}_k+\hat{\bm{\eta}}_k:=\hat{C}\hat{\bm{\Psi}}_k+\hat{\bm{\eta}}_k
  \end{align} 
  \end{subequations}
 where $\hat{\bm{\Psi}}_k=[
               \hat{\bm{\psi}}_k~\hat{\dot{\bm{\psi}}}_k
            ]^T$ and $\hat{D}_{kk}=-4\sin^2\left(\frac{\pi k}{n}\right)$.
           The solution to  (\ref{eq:KFfourierss}, \ref{eq:nondimKFcostfourier}) is given by
           $\dot{\hat{\bm{\Psi}}}_k=(\hat{A}_k-\hat{L}_k\hat{C})\hat{\bm{\Psi}}_k+\hat{L}_k\hat{\bm{\gamma}}_k$ where 
           \begin{align}
               \hat{L}_k=\hat{S}_k\hat{C}^T\hat{Q}_k &=\begin{bmatrix}
                  \hat{S}_1(k)&\hat{S}_0(k)\\\hat{S}_0(k)&\hat{S}_2(k)
                  \end{bmatrix}\begin{bmatrix}
                   \Pi_4\notag\\0
               \end{bmatrix}(1-\Pi_1\hat{D}_{kk})\label{eq:genL}\\
            \end{align}
            and $\hat{S}_k \succ 0$ is the solution to the Algebraic Riccati Equation \begin{equation}\label{eq:KFARE}
               \hat{A}_k\hat{S}_k+\hat{S}_k\hat{A}^T_k+\hat{B}\hat{B}^T-\hat{S}_k\hat{C}^T\hat{Q}_k\hat{C}\hat{S}_k=0.
            \end{equation}
         Solving (\ref{eq:KFARE}) gives 
         \begin{align}\hat{S}_0(k)&=\frac{\hat{D}_{kk}\pm\sqrt{\hat{D}^2_{kk}+\Pi_4^2\left(1-\Pi_1\hat{D}_{kk}\right)}}{\Pi_4^2(1-\Pi_1\hat{D}_{kk})}\\
         \hat{S}_1(k)&=\sqrt{\frac{2\hat{S}_0(k)}{{\Pi_4^2(1-\Pi_1\hat{D}_{kk})}}}\end{align} so that (\ref{eq:genL}) can be written as (\ref{eq:finalL}).
\end{proof}
\end{document}